\documentclass{gtart_h}

\def\ifplaintex{\expandafter\ifx\csname documentclass\endcsname\relax}

\def\gtp{{\mathsurround=0pt\it $\cal G\mskip-2mu$eometry \&\ 
$\cal T\!\!$opology $\cal P\!$ublications}}  

\def\recd{{\small Received:\qua\receiveddate\ifx\reviseddate\relax
\else\qquad Revised:\qua\reviseddate\fi\par}} 

\def\lognumber#1{\def\thelognumber{#1}}
\def\volumenumber#1{\def\thevolumenumber{#1}}
\def\volumeyear#1{\def\thevolumeyear{#1}}
\def\papernumber#1{\def\thepapernumber{#1}}
\def\pagenumbers#1#2{\def\startpage{#1}\def\finishpage{#2}}
\def\published#1{\def\publishdate{#1}}

\def\received#1{\def\receiveddate{#1}}

\def\accepted#1{\def\accepteddate{#1}}
\def\asciititle#1{\def\theasciititle{#1}}

\def\asciiauthors#1{\def\theasciiauthors{#1}}
\def\asciiaddress#1{\def\theasciiaddress{#1}}
\def\asciiemail#1{\def\theasciiemail{#1}}

\def\coverauthors#1{\def\thecoverauthors{#1}}
\long\def\asciiabstract#1{\long\def\theasciiabstract{#1}}
\def\asciikeywords#1{\def\theasciikeywords{#1}}

\let\\\par\let\thelognumber\relax\let\thevolumenumber\relax
\let\thepapernumber\relax\let\thevolumeyear\relax\let\startpage\relax
\let\finishpage\relax\let\publishdate\relax\let\receiveddate\relax
\let\reviseddate\relax\let\accepteddate\relax\let\theasciititle\relax
\let\theasciiauthors\relax\let\theasciiaddress\relax
\let\theasciiabstract\relax\let\theasciikeywords\relax

\let\thecoverauthors\relax\let\theasciiemail\relax

\ifplaintex
\font\logobig=cmssbx10 scaled 3836
\font\logomed=cmssbx10 scaled 2557
\else
\font\logobig=cmssbx10 scaled 4200
\font\logomed=cmssbx10 scaled 2800
\fi

\long\def\makeagttitle{   
\count0=\startpage
\agt\hfill      
\hbox to 45truept{\vbox to 0pt{\vglue -13truept{\logomed A\kern -.37em{\logobig 
T}\kern -.38em G}\vss}\hss}
\break
{\small Volume \thevolumenumber\ (\thevolumeyear)
\startpage--\finishpage\nl
Published: \publishdate}

\vglue .25truein

{\parskip=0pt\leftskip 0pt plus
1fil\def\\{\par\smallskip}{\Large\bf\thetitle}\par\medskip} \vglue
0.05truein

{\parskip=0pt\leftskip 0pt plus 1fil\def\\{\par}{\sc\theauthors}
\par\medskip}%
 
\vglue 0.03truein 

{\small\leftskip 25truept\rightskip 25truept{\bf Abstract}\stdspace\theabstract

{\bf AMS Classification}\stdspace\theprimaryclass
\ifx\thesecondaryclass\relax\else; \thesecondaryclass\fi\par
{\bf Keywords}\stdspace \thekeywords\par}\vglue 7truept

}   

\ifplaintex
\font\phead=cmsl9 scaled 950
\font\pnum=cmbx10 scaled 913
\font\pfoot=cmsl9 scaled 950
\headline{\vbox to 0pt{\vskip -4.5mm\line{\small\phead\ifnum
\count0=\startpage ISSN 1472-2739 (on-line) 1472-2747 (printed)
\hfill {\pnum\folio}\else\ifodd\count0\def\\{ }%
\ifx\theshorttitle\relax\thetitle\else\theshorttitle\fi\hfill{\pnum\folio}
\else\def\\{ and }{\pnum\folio}\hfill\ifx\theshortauthors\relax\theauthors
\else\theshortauthors\fi\fi\fi}\vss}}
\footline{\vbox to 0pt{\vglue 0mm\line{\small\pfoot\ifnum\count0=\startpage
\copyright\ \gtp\hfill\else
\agt, Volume \thevolumenumber\ (\thevolumeyear)\hfill\fi}\vss}}
\else
\font=cmsl9 scaled 1050
\font\lnum=cmbx10 
\font=cmsl9 scaled 1050
\makeatletter
\def\@oddhead{{\small\lhead\ifnum\count0=\startpage ISSN 1472-2739 
(on-line) 1472-2747 (printed)\hfill {\lnum\number\count0}\else\ifodd\count0
\def\\{ }\ifx\theshorttitle\relax \thetitle \else\theshorttitle\fi\hfill
{\lnum\number\count0}\else\def\\{ and }{\lnum\number\count0}
\hfill\ifx\theshortauthors\relax 
\theauthors\else\theshortauthors\fi\fi\fi}}\def\@evenhead{\@oddhead}
\def\@oddfoot{\small\lfoot\ifnum\count0=\startpage\copyright\ \gtp\hfill\else
\agt, Volume \thevolumenumber\ (\thevolumeyear)\hfill\fi}
\def\@evenfoot{\@oddfoot}
\makeatother
\fi
\let\maketitlepage\makeagttitle
\let\makeshorttitle\maketitlepage
\let\maketitle\maketitlepage

\newwrite\gtoutfile
\long\gdef\makeheadfile{  
{\def\\{, }\def\s{ }
\immediate\openout\gtoutfile head.xxx
\immediate\write\gtoutfile{Proxy-for: \ifx\theasciiauthors\relax
\theauthors\else\theasciiauthors\fi\s<\ifx\theasciiemail\relax\theemail\else\theasciiemail\fi>}
\immediate\write\gtoutfile{\noexpand\\}
\immediate\write\gtoutfile{Authors: \ifx\theasciiauthors\relax
\theauthors\else\theasciiauthors\fi}
{\def\\{ }\immediate\write\gtoutfile{Title: \ifx\theasciititle\relax
\thetitle\else\theasciititle\fi}}
\immediate\write\gtoutfile{Subj-class: GT or SG, GR etc}
\immediate\write\gtoutfile{MSC-class: \theprimaryclass\ifx\thesecondaryclass\relax\else, \thesecondaryclass\fi}
\immediate\write\gtoutfile{Journal-ref: Algebr. Geom. Topol. \thevolumenumber\s
(\thevolumeyear) \startpage-\finishpage}
\immediate\write\gtoutfile{Comments: Published by Algebraic and
Geometric Topology at}
\immediate\write\gtoutfile{\s\s\s  http://www.maths.warwick.ac.uk/agt/AGTVol\thevolumenumber/agt-\thevolumenumber-\thepapernumber.abs.html}
\immediate\write\gtoutfile{\noexpand\\}
\immediate\write\gtoutfile{}
\ifx\theasciiabstract\relax
\immediate\write\gtoutfile{\theabstract}\else
\immediate\write\gtoutfile{\theasciiabstract}\fi
\immediate\write\gtoutfile{}
\immediate\write\gtoutfile{\noexpand\\}
\immediate\write\gtoutfile{}
\immediate\closeout\gtoutfile}}  

\def\maketitlepage{\makeagttitle\makeheadfile}
\let\makeshorttitle\maketitlepage
\let\maketitle\maketitlepage

\lognumber{49}
\volumenumber{4}
\volumeyear{2004}
\papernumber{49}
\pagenumbers{1125}{1144}
\received{21 June 2004} 
\accepted{4 November 2004}
\published{25 November 2004}

\usepackage{amssymb, amsmath, graphs}

\newtheorem{thm}{Theorem}[section]
\newtheorem{cor}[thm]{Corollary}
\newtheorem{lem}[thm]{Lemma}
\newtheorem{prop}[thm]{Proposition}
\theoremstyle{definition}
\newtheorem{defi}[thm]{Definition}

\numberwithin{equation}{section}
\newcommand{\al}{\alpha}
\newcommand{\be}{\beta}
\newcommand{\ga}{\gamma}
\newcommand{\Ga}{\Gamma}

\newcommand{\de}{\delta}

\newcommand{\Om}{\Omega}
\newcommand{\om}{\omega}

\newcommand{\si}{\sigma}
\newcommand{\Si}{\Sigma}
\renewcommand{\th}{\theta}

\newcommand{\s}{\mathbf s}

\newcommand{\FF}{\mathcal F}

\newcommand{\Z}{\mathbb Z}
\newcommand{\Q}{\mathbb Q}

\newcommand{\CP}{{\mathbb C}{\mathbb P}}

\newcommand{\del}{\partial}

\newcommand{\x}{\times}

\begin{document}

\title{Transverse contact structures\\on Seifert 3--manifolds}
\asciititle{Transverse contact structures\\on Seifert 3-manifolds}

\author{Paolo Lisca\\Gordana Mati\'c}
\coverauthors{Paolo Lisca\\Gordana Mati\noexpand\'c}
\asciiauthors{Paolo Lisca\\Gordana Matic}

\address{Dipartimento di Matematica ``L.~Tonelli''\\ 
Universit\`a di Pisa \\I-56127
Pisa, ITALY\\{\rm and}\\Department of Mathematics\\ University of Georgia\\Athens, GA
30602, USA}
\asciiaddress{Dipartimento di Matematica "L. Tonelli"\\ 
Universita di Pisa \\I-56127
Pisa, ITALY\\and\\Department of Mathematics\\University of Georgia\\Athens, GA
30602, USA}

\gtemail{\mailto{lisca@dm.unipi.it}, \mailto{gordana@math.uga.edu}}
\asciiemail{lisca@dm.unipi.it, gordana@math.uga.edu}

\begin{abstract} 
We characterize the oriented Seifert--fibered three--manifolds which 
admit positive, transverse contact structures.
\end{abstract}
\asciiabstract{%
We characterize the oriented Seifert-fibered three-manifolds which 
admit positive, transverse contact structures.}
\asciikeywords{Transverse contact structures, Seifert three-manifolds}

\maketitle

\section{Introduction and statement of results}
\label{s:intro}

Foliations and contact structures are arguably at opposite ends of the
spectrum of the possible 2-plane fields $\xi$ on a 3-manifold. While
foliations are integrable plane fields, contact structures are totally
non-integrable. If the planes in the distribution $\xi$ are given as
kernels of a one form $\alpha$, then they form a foliation if $\alpha
\wedge d\alpha = 0$, while they form a contact structure if everywhere
pointwise $\alpha \wedge d\alpha \neq 0$.  Eliashberg and
Thurston~\cite{ET} introduced the notion of confoliation, which
generilizes the notion of foliation as well as the notion of contact
structure. They also showed that a taut foliation can be
$C^0$--approximated by a fillable contact structure.

The problem of existence of foliations transverse to the $S^1$--action
on a Seifert fibered 3-manifold $M$ was studied, and almost completely
solved, in the 1980's by Eisenbud, Hirsch, Jankins and Neumann who
translated, via holonomy, the original problem into a problem about
homeomorphisms of the circle~\cite{EHN, JN1, JN2}. The problem was
settled by Naimi~\cite{Na}.

The existence of transverse contact structures is clearly a related
question.  A foliation transverse to the $S^1$--action on a Seifert
3--manifold is taut since any orbit provides a closed loop transverse
to the leaves of the foliation. If a transverse foliation $\FF$ exists
and the underlying manifold is not $S^1\x S^2$ then, according to
Eliashberg and Thurston~\cite{ET}, $\FF$ can be approximated by a
fillable contact structure $\xi$. The approximation is still
transverse to the circle action on $M$.  However, often $M$ supports a
transverse tight contact structure even when there are no transverse
taut foliations.

Let $p\co Y\to\Si_g$ be an oriented three--dimensional circle bundle
with Euler class $e(Y)$ and base of genus $g$. Here $g\in\Z$ is
defined so that
\[
\chi(\Si_g)=
\begin{cases}
2-2g\quad\text{if $\Si_g$ is orientable},\\
2+g\quad\text{if $\Si_g$ is non--orientable},
\end{cases}
\]
where $\chi(\Si_g)$ is the Euler characteristic of $\Si_g$. Note that,
with this convention, $\Si_g$ is orientable for $g\geq 0$ and
non--orientable for $g<0$. Milnor and Wood studied the case of
foliations transverse to the fibers of $p$. According to~\cite{Mi,
Wo}, $Y$ carries a transverse foliation if and only if one of the
following holds:
\begin{itemize}
\item
$\chi(\Si_g)\leq 0$ and  $|e(Y)|\leq -\chi(\Si_g)$,
\item
$\chi(\Si_g)\geq 0$ and $e(Y)=0$.
\end{itemize}
More recently, Giroux~\cite{Gi} and Sato and Tsuboi~\cite{ST} proved
that $Y$ admits a positive, transverse contact structure if and only
if one of the following holds:
\begin{itemize}
\item
$\chi(\Si_g)\leq 0$ and $e(Y)\leq -\chi(\Si_g)$,
\item
$\chi(\Si_g)>0$ and $e(Y) < 0$.
\end{itemize}
Let $p\co M\to\Si_g$ be an oriented three--dimensional Seifert
fibration with base of genus $g$ and normalized (in the sense
of~\cite{Or}) Seifert invariants
\begin{equation*}
\{b, g; (\al_1,\be_1),\ldots, (\al_r,\be_r)\}. 
\end{equation*}
Define 
\[
e(M):=-b-\sum_{i=1}^r \frac{\be_i}{\al_i},\quad
e_0(M):=-b-r,\quad\text{and}\quad\Ga(M):=(\ga_1,\ldots,\ga_r),
\]
where 
\[
\ga_i:=1-\frac{\be_i}{\al_i},\quad i=1,\dots,r.
\]
Recall that the normalized Seifert invariants of $-M$ are 
\[
\{-b-r, g; (\al_1,\al_1-\be_1),\ldots, (\al_r,\al_r-\be_r)\}.
\]
Therefore, 
\[
e(-M)=-e(M),\quad e_0(-M)=-e_0(M)-r=b,\quad
\Ga(-M)=(1-\ga_1,\ldots,1-\ga_r).
\]
\begin{defi}
We say that $\Ga=(\ga_1,\ldots,\ga_r)\in(\Q\cap(0,1))^r$
is~\emph{realizable} if $r\geq 3$ and there exist a permutation
$\si\in{\cal S}_r$ and coprime integers $m>a>0$ such that:
\[
\ga_{\si(1)} < \frac am, \quad
\ga_{\si(2)} < \frac{m-a}m, \quad
\ga_{\si(3)},\ldots,\ga_{\si(r)} < \frac 1m.
\]
\end{defi}
The main results of Eisenbud, Hirsch, Jankins, Neumann and Naimi on
the existence of transverse foliations can be summarized in the
following statement.

\begin{thm}{\rm\cite{EHN, JN2, Na}}\label{t:quoted}\qua
Let $p\co M\to\Si_g$ be an oriented three--dimensional Seifert
fibration as above. Then, $M$ carries a smooth foliation transverse to
the fibration if and only if one of the following holds:
\begin{enumerate}
\item[\rm (a)]
$e_0(M)\leq -\chi(\Si_g)$ and $e_0(-M)\leq -\chi(\Si_g)$
\item[\rm (b)]
$g=0$ and $e(M)=e(-M)=0$ 
\item[\rm (c)]
$g=0$, $e_0(M)=-1$ and $\Ga(M)$ is realizable
\item[\rm (d)]
$g=0$, $e_0(-M)=-1$ and $\Ga(-M)$ is realizable
\end{enumerate}
\end{thm}

\rk{\bf Warning} The reader should be aware of a difference in
conventions: the quantity denoted by `$b$' in~\cite{JN2} corresponds
to $-b$ in our present notation.

The main result of this paper is the following theorem, giving
necessary and sufficient conditions for the existence of contact
structures transverse to the fibers of a Seifert fibration:

\begin{thm}\label{t:main}
Let $p\co M\to\Si_g$ be an oriented three--dimensional Seifert
fibration as above. Then, $M$ carries a positive contact structure
transverse to the fibers of $p$ if and only if one of the following
holds:
\begin{enumerate}
\item[\rm (a)]
$e_0(M)\leq -\chi(\Si_g)$
\item[\rm (b)]
$g=0$, $r\leq 2$ and $e(M)<0$
\item[\rm (c)]
$g=0$, $e_0(M)=-1$ and $\Ga(M)$ is realizable.
\end{enumerate}
\end{thm}

In view of the Eliashberg--Thurston's approximation theorem, one might
wonder why Condition (b) from Theorem~\ref{t:quoted} does not appear
explicitely in the statement of Theorem~\ref{t:main}. The reason is
that $e(M)=e(-M)=0$ implies $e_0(M)\leq -1$ and $e_0(-M)\leq -1$. If
both inequalities are strict, case (a) of Theorem~\ref{t:main} holds
for both $M$ and $-M$. If both inequalities are equalities, then $r=2$
and it is easy to check that $M=S^1\x S^2$, therefore the
Eliashberg--Thurston's perturbation theorem cannot be applied. If
$e_0(M)=-1$ and $e_0(-M)\leq -2$, (a) of Theorem~\ref{t:main} holds
for $-M$ and, since $r\geq 3$, it follows from~\cite{JN2, Na} that (c)
of Theorem~\ref{t:main} holds for $M$.  Similarly, if $e_0(M)=-2$ and
$e_0(-M)\leq -1$ then (a) of Theorem~\ref{t:main} holds for $M$ and
(c) of Theorem~\ref{t:main} for $-M$.

It should be mentioned that Honda~\cite{H} uses monodromy arguments
(parallel to the ones used in the foliation case) to study the
existence of transverse contact structures, obtaining results similar
to ours.

Our approach is completely different. We show that the conditions of
Theorem~\ref{t:main} are imposed on a Seifert 3--manifold $M$ by the
existence of a transverse contact structure $\xi$ because of a certain
property of a specific symplectic filling $X_M$ of $(M,\xi)$. The
property is simply the fact that surfaces embedded in $X_M$ and any of
its blowdowns satisfy the adjunction inequalities.

To prove Theorem~\ref{t:main} when $\Si_g$ is orientable, we look at
$M$ as the boundary of an equivariant plumbing $X_M$ built according
to a weighted star $\Gamma$ of Figure~\ref{f:fig0} as described
in~\cite{Or}. McCarthy and Wolfson~\cite{McW} showed that $X_M$
carries a symplectic structure such that the orbits of the
$S^1$--action on $M$ are tangent to the kernel of the restriction of
the symplectic form. For any transverse contact structure $\xi$ on
$M$, this makes $X_M$ a symplectic filling of $(M,\xi)$. To prove that
the conditions of Theorem~\ref{t:main} are necessary, we impose that
the adjunction inequalities be satisfied by surfaces embedded in $X_M$
and, if the central node of $\Gamma$ is a $(-1)$-sphere, into
successive blow-downs of $X_M$.

It turns out that the most difficult case to analyze is case (c) of
Theorem~\ref{t:main}, involving the ``realizability" condition. It is
very interesting to see it appear from this new angle. In fact, one
may think of the proof of Theorem~\ref{t:main} as a way of checking
whether a given $\Ga(M)=(\ga_1,\ldots,\ga_r)$ is realizable. In short,
$\Ga(M)$ is realizable if and only if in the successive blow-downs of
the plumbed manifold $X_M$ there are no obvious surfaces violating the
adjunction inequality.

The case when $\Si_g$ is non--orientable is reduced to the orientable
case by pulling back the fibration $M\to\Si_g$ to the orientable 
double cover of $\Si_g$.

To prove the existence of a transverse contact structure under the
assumptions of Theorem~\ref{t:main}, we either use the results on
foliations and the Eliashberg-Thurston approximation result, or we
construct the contact structure directly.

The paper is organized as follows. In Section~\ref{s:filling} we
describe in detail the symplectic filling $X_M$, and we collect its
necessary properties. In Section~\ref{s:invariant} we characterize
those Seifert 3--manifolds which admit $S^1$--invariant transverse
contact structures, and in Section~\ref{s:mainproof} we prove
Theorem~\ref{t:main}.

\section{Symplectic fillings}
\label{s:filling}

Let $p\co M\to\Si_g$ be a Seifert $3$--manifold  with
orientable base and Seifert invariants
\[
\{b,g;(\al_1,\be_1),\ldots,(\al_r,\be_r)\}.
\]
Suppose that 
\[
\frac{\al_i}{\al_i-\be_i} = [b_{i,1},\ldots, b_{i,s_i}], \quad
i=1,\ldots, r,
\]
where 
\[
[b_1,\ldots, b_k] := 
b_1 - \cfrac{1}{b_2 -
          \cfrac{1}{\ddots -
             \cfrac{1}{b_k
}}},\quad
b_1,\ldots,b_k\geq 2
\]
Note that, since we assume $0 < \be_i< \al_i$, there is a unique
continued expansion with all $b_{i,j} \geq 2$.

\begin{figure}[hb]
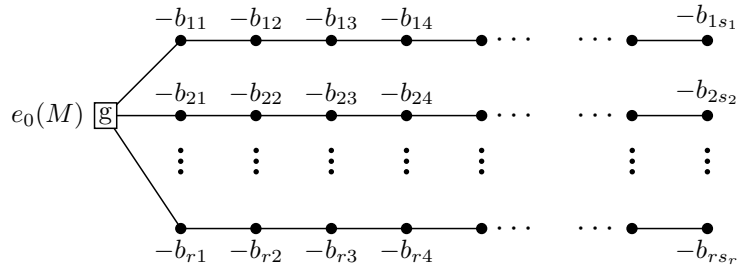
\small
\begin{center}
\begin{graph}(10,3.5)(-9,-2)
\graphnodesize{0.15}
  \textnode{n0}(-8,0){g}
  \roundnode{n11}(-7,1)
  \roundnode{n12}(-6,1)
  \roundnode{n13}(-5,1)
  \roundnode{n14}(-4,1)
  \roundnode{n15}(-3,1)
  
  \roundnode{n16}(-1,1)
  \roundnode{n17}(0,1)

  \roundnode{n21}(-7,0)
  \roundnode{n22}(-6,0)
  \roundnode{n23}(-5,0)
  \roundnode{n24}(-4,0)
  \roundnode{n25}(-3,0)
  
  \roundnode{n26}(-1,0)
  \roundnode{n27}(0,0)

  \roundnode{nr1}(-7,-1.5)
  \roundnode{nr2}(-6,-1.5)
  \roundnode{nr3}(-5,-1.5)
  \roundnode{nr4}(-4,-1.5)
  \roundnode{nr5}(-3,-1.5)
  
  \roundnode{nr6}(-1,-1.5)
  \roundnode{nr7}(0,-1.5)

  \edge{n0}{n11}
  \edge{n11}{n12}
  \edge{n12}{n13}
  \edge{n13}{n14}
  \edge{n14}{n15}
  \edge{n16}{n17}
 
  \edge{n0}{n21}
  \edge{n21}{n22}
  \edge{n22}{n23}
  \edge{n23}{n24}
  \edge{n24}{n25}
  \edge{n26}{n27}

  \edge{n0}{nr1}
  \edge{nr1}{nr2}
  \edge{nr2}{nr3}
  \edge{nr3}{nr4}
  \edge{nr4}{nr5}
  \edge{nr6}{nr7}

  \autonodetext{n0}[w]{\small $e_0(M)$}
  \autonodetext{n11}[n]{$-b_{11}$}
  \autonodetext{n12}[n]{$-b_{12}$}
  \autonodetext{n13}[n]{$-b_{13}$}
  \autonodetext{n14}[n]{$-b_{14}$}
  \autonodetext{n15}[e]{\large $\cdots$}
  \autonodetext{n16}[w]{\large $\cdots$}
  \autonodetext{n17}[n]{$-b_{1s_1}$}

  \autonodetext{n21}[n]{$-b_{21}$}
  \autonodetext{n22}[n]{$-b_{22}$}
  \autonodetext{n23}[n]{$-b_{23}$}
  \autonodetext{n24}[n]{$-b_{24}$}
  \autonodetext{n25}[e]{\large $\cdots$}
  \autonodetext{n26}[w]{\large $\cdots$}
  \autonodetext{n27}[n]{$-b_{2s_2}$}

  \autonodetext{nr1}[s]{$-b_{r1}$}
  \autonodetext{nr2}[s]{$-b_{r2}$}
  \autonodetext{nr3}[s]{$-b_{r3}$}
  \autonodetext{nr4}[s]{$-b_{r4}$}
  \autonodetext{nr5}[e]{\large $\cdots$}
  \autonodetext{nr6}[w]{\large $\cdots$}
  \autonodetext{nr7}[s]{$-b_{rs_r}$}

  \freetext(-7,-0.5){\Huge $\vdots$}
  \freetext(-6,-0.5){\Huge $\vdots$}
  \freetext(-5,-0.5){\Huge $\vdots$}
  \freetext(-4,-0.5){\Huge $\vdots$}
  \freetext(-3,-0.5){\Huge $\vdots$}
  \freetext(-1,-0.5){\Huge $\vdots$}
  \freetext(0,-0.5){\Huge $\vdots$}
\end{graph}
\end{center}
\caption{The plumbing $X_M$}
\label{f:fig0} 
\end{figure}

By~\cite[Chapter~2]{Or}, $M$ is isomorphic (as a $3$--manifold with
$S^1$--action) to the boundary of a $4$--manifold with $S^1$--action
$X_M$ obtained by equivariant plumbing according to the graph $\Ga$ of
Figure~\ref{f:fig0}, where the central vertex represents a disk bundle
with Euler number $e_0(M)$ over a surface of genus $g$, and all the
other vertices represent disk bundles over the sphere with Euler
numbers $-b_{i,j}$.

\begin{thm}[\cite{McW}]\label{t:McW} 
The $S^1$--manifold $X_M$ carries a symplectic form $\om$ such that
every orbit of the $S^1$--action on $\del X_M$ is tangent to the
kernel of $\om|_{\del X_M}$.\qed
\end{thm}

Recall that a {\sl symplectic filling} of a closed contact
3--manifold $(M,\xi)$ is a symplectic 4--manifold $(X,\om)$ such that
(i) $X$ is oriented by $\om\wedge\om$, (ii) $\del X=M$ as oriented
manifolds, and (iii) $\om|_\xi\neq 0$ at every point of $M$.

\begin{cor}
\label{c:transverse-fillable}
Let $p\co M\to\Si_g$ be an oriented Seifert fibered $3$--manifold with
orientable base $\Si_g$. Let $\xi$ be a a positive contact structure
on $M$ transverse to the fibration. Then, $(X_M,\om)$ from
Theorem~\ref{t:McW} is a symplectic filling of $(M,\xi)$. Moreover,
$\xi$ is universally tight.
\end{cor}

\begin{proof}
The first part follows immediately from the definitions together with
Theorem~\ref{t:McW}.  For the second part, recall that the fundamental
group of a Seifert $3$--manifold is residually
finite~\cite[pp.~176--177]{He}. Therefore, in order to prove that the
contact structure is universally tight it is enough to show that the
pull--back of the contact structure to every finite cover is
tight. This follows from the first part of the statement, because
every finite cover is still Seifert fibered, the pulled--back structure
is transverse, and fillable contact structures are tight~\cite{Eli1}.
\end{proof}

\begin{cor}\label{c:adjunction}
Let $p\co M\to\Si_g$ be an oriented Seifert fibered $3$--manifold with
orientable base $\Si_g$. Suppose that $M$ carries a positive contact
structure transverse to the fibration. Let $S\subset X_M$ be an
oriented surface of genus $g(S)$ smoothly embedded in the 4--manfild
$X_M$ given by Figure~\ref{f:fig0}. Then,
\[
\begin{cases}
S\cdot S \leq -1\quad\text{if}\quad g(S)=0,\\
S\cdot S \leq -\chi(S)\quad\text{if}\quad g(S)>0.
\end{cases}
\]
\end{cor}

\begin{proof}
Let $\xi$ be a a positive contact structure on $M$ transverse to the
fibration. By Corollary~\ref{c:transverse-fillable}, $(M,\xi)$ has a
symplectic filling of the form $(X_M,\om)$. By~\cite{El}, $(X_M,\om)$
can be compactified to a closed, symplectic 4--manifold $\widehat
X_M$.  Moreover, up to adding a suitable Stein cobordism to $X_M$
(using e.g.~\cite[Theorem~2.5]{EH}) it can be arranged that
$b^+_2(\widehat X_M)>1$.  Since the stated inequalities are
satisfied by surfaces in $\widehat X_M$~\cite{OS}, the conclusion
follows.
\end{proof}

\section{$S^1$--invariant structures}
\label{s:invariant}

This section consists of the following proposition, which
characterizes the oriented Seifert 3--manifolds carrying positive,
$S^1$--invariant transverse contact structures.

\begin{prop}\label{p:e<0}
Let $p\co M\to\Si_g$ be an oriented Seifert 3--manifold. Then, $M$
carries a positive, $S^1$--invariant transverse contact structure if
and only if $e(M)<0$.
\end{prop}

\begin{proof}
Suppose that $M$ has Seifert invariants
\[
\{b, g; (\al_1,\be_1),\ldots, (\al_r,\be_r)\}.
\]
Let $a=\al_1\cdots\al_r$. The cyclic group $\Z/a\Z$ acts on $M$ when
it is regarded as a subgroup of $S^1$. The quotient $M'=M/\Z/a\Z$ is a
genuine $S^1$--bundle over the same base $\Si_g$. Moreover,
by~\cite[Theorem~1.2]{NR},
\begin{equation}\label{e:cover}
e(M')=a e(M).
\end{equation}
Thus, if $e(M)<0$ then $e(M')<0$ and by~\cite{Gi, ST} $M'$ carries an
$S^1$--invariant contact structure $\xi$ transverse to the fibers. The 
pull--back of $\xi$ to $M$ is also $S^1$--invariant and tranverse, 
therefore we have proved the first half of the statement. 

Now let $\xi$ be an $S^1$--invariant transverse contact structure on
$M$.  We claim that there exists an induced $S^1$--invariant
transverse contact structure $\overline\xi$ on $M'$. By~\cite{Gi, ST}
and Equation~\eqref{e:cover}, proving the claim clearly suffices to
finish the proof.

To prove the claim, we first argue locally around a singular
fiber. Recall that a neighborhood of a singular fiber $F$ is of the
form
\[
(D^2\x S^1)/\Z/p\Z,
\]
where the generator $g\in\Z/p\Z$ acts by $(x,t)\mapsto (xg,tg^q)$ for
some $0<q<p$ coprime with $p$. The contact structure $\xi$ lifts to a
contact structure $\tilde\xi$ on $D\x S^1$ which is
$\Z/p\Z$--invariant and $S^1$--invariant under the standard
$S^1$--action on the second factor.  Since $\xi$ is transverse to the
singular fiber $F$, after a suitable change of trivialization of the
neighborhood of $F$ $\tilde\xi$ is given as the zero set of a
$1$--form $\al+d\th$, where $d\th$ is the standard angular $1$--form
on $S^1$, $\al\in\Om^1(D^2)$ is $S^1$--invariant and
$d\al\in\Om^2(D^2)$ is a volume form.

Since the $\Z/p\Z$--action and the $\Z/a\Z$--action commute, the quotient
map from $D^2\x S^1$ onto a neighorhood of the image $F'\subset M'$ of
$F$ can be factored as follows:
\[
D^2\x S^1\stackrel{f}{\to}(D^2\x S^1)/\Z/a\Z\stackrel{g}{\to}
\left((D^2\x S^1)/\Z/a\Z\right)/\Z/p\Z.
\]
Here $\Z/a\Z$ acts freely on the second factor as a subgroup of $S^1$,
and after the identification $(D^2\x S^1)/\Z/a\Z\cong D^2\x S^1$,
$\Z/p\Z$ acts by rotations on the first factor. Since $\al$ and $d\th$
are $S^1$--invariant, $\tilde\xi$, and therefore $\xi$, descends to a
smooth contact structure $\overline\xi$ on a neighborhood of
$F'$. Since $\xi$ also descends to a transverse $S^1$--invariant
contact structure when restricted to the complement of the the
singular fibers, the claim is proved.
\end{proof}
\eject

\section{The proof of Theorem~\ref{t:main}}\label{s:mainproof}

\subsection{Sufficiency of the conditions}

\sh{The case $\chi(\Si_g)\leq 0$}

Suppose that the inequality stated in case (a) of Theorem~\ref{t:main}
holds. Let us first assume that $e_0(-M)\leq -\chi(\Si_g)$. Then, it
follows from case (a) of Theorem~\ref{t:quoted} that $M$ carries a
smooth foliation $\FF$ transverse to the Seifert fibration. Moreover,
$M\neq S^1\x S^2$ because $\chi(\Si_g)\leq 0$ implies that the
fundamental group of $\Si_g$ is non--trivial, and this is incompatible
with $M=S^1\x S^2$. Therefore, by~\cite{ET} there exist positive
contact structures arbitrarily $C^0$--near $\FF$, hence transverse to
the fibration.

If $e_0(-M)>-\chi(\Si_g)$, take $N$ to be the Seifert fibered
$3$--manifold with Seifert invariants 
\[
(-\chi(\Si_g), g;(\al_1,\be_1),\ldots, (\al_r,\be_r)).  
\]
Since 
\[
e_0(N)=\chi(\Si_g)-r\leq -\chi(\Si_g), 
\]
$N$ carries a positive contact structure transverse to the
fibration. Moreover, $M$ can be obtained from $N$ by a smooth surgery
along a regular fiber $f\subset N$ with coefficient
\[
\frac{1}{e_0(-M)+\chi(\Si_g)}
\]
with respect to the framing induced by the fibration. It is not
difficult to check(see~\cite[Lemma~1.3]{Gi}) that transverse contact
structures extend from the complement of a regular neighborhood
of $f$ in $N$ to transverse contact structures on $M$.

\sh{The case $\chi(\Si_g)>0$, i.e.~$g\in\{-1,0\}$}

If case (a) of Theorem~\ref{t:main} holds, then
\begin{equation}\label{e:inequa}
e(M) = e_0(M) + \sum_{i=1}^r (1-\frac{\be_i}{\al_i}) \leq -\chi(\Si_g)
+ \sum_{i=1}^r (1-\frac{\be_i}{\al_i}).
\end{equation}
If $e(M)<0$, then $M$ carries a transverse contact structure by
Proposition~\ref{p:e<0}. If $e(M)\geq 0$, then
Equation~\eqref{e:inequa} implies $r\geq 2$ if $g=-1$, and $r\geq 3$
if $g=0$, therefore $M\neq S^1\x S^2$. Also, since $e(-M)=-e(M)\leq
0$,
\[
e_0(-M) = e(-M) - \sum_{i=1}^r\frac{\be_i}{\al_i} < 0,
\]
therefore $e_0(-M)\leq -1$. If (i) $g=-1$ or (ii) $g=0$ and
$e_0(-M)\leq -2$, then by Theorem~\ref{t:quoted}(a) together with the
Eliashberg and Thurston's theorem we are done. If $g=0$ and
$e_0(-M)=-1$ then, since $e(-M)\leq 0$, by~\cite[Theorem~1]{JN1} (see
also~\cite[Theorem~3.1]{JN2}) and~\cite[Theorem~3]{JN2}, $-M$, and
therefore $M$, supports a smooth foliation $\FF$ transverse to the
fibration and we conclude as before.

If case (b) of Theorem~\ref{t:main} holds, then the conclusion follows
by Proposition~\ref{p:e<0}.

If (c) holds, then we get a transverse contact structure via
Theorem~\ref{t:quoted}(c) and~\cite{ET}.

\subsection{Necessity of the conditions}\label{s:necessity}

\sh{The case $\chi(\Si_g)\leq 0$}

If $g>0$, observe that the smooth 4--manifold $X_M$ of
Section~\ref{s:filling} contains a smooth surface $\Si_g$ of genus $g$
and self--intersection $e_0(M)$. If $M$ carries a positive, transverse
contact structure, then by Corollary~\ref{c:adjunction} $e_0(M)\leq
-\chi(\Si_g)$. Therefore case (a) of Theorem~\ref{t:main} holds.

If $g<-1$, let $\widehat\Si_g\to\Si_g$ be the orientable double cover
of $\Si_g$, and let $\widehat M\to\widehat\Si_g$ be the pull--back of
the fibration $M\to\Si_g$. According to~\cite[Theorem~1.2]{NR},
$e_0(\widehat M)=2e_0(M)$. If $M$ carries a positive, transverse
contact structure then so does $\widehat M$. Therefore, since
$\widehat\Si_g$ has positive genus we have
\[
e_0(\widehat M)\leq -\chi(\widehat\Si_g)=-2\chi(\Si_g),
\]
hence case (a) of Theorem~\ref{t:main} holds.

\sh{The case $\chi(\Si_g)>0$}

Suppose first that $g=0$. Then, Corollary~\ref{c:adjunction} implies
$e_0(M)\leq -1$. Hence, either case (a) of Theorem~\ref{t:main} holds, 
or $e_0(M)=-1$.  Let us assume the latter.

If $r=0$ then $e(M)=e_0(M)=-1$, therefore case (b) of
Theorem~\ref{t:main} holds. If $r=1$ or $r=2$, we blow down
$(-1)$--spheres in $X_M$ as far as possible. The adjuction
inequalities imply that no non--negative sphere appears, therefore
$X_M$ is negative definite and, since $e(M)$ is an eigenvalue of
$Q_X$~\cite[Theorem~5.2]{NR}, $e(M)<0$. Thus, again case (b) of
Theorem~\ref{t:main} holds.

When $g=0$, we are left to consider the subcase $g=0$ and $r\geq
3$. Before tackling this subcase we deal with the case $g=-1$.

If $g=-1$, let $\widehat M\to\widehat\Si_g$ be the pull--back of
$M\to\Si_g$ under the orientable double cover of
$\Si_g$. By~\cite[Theorem~1.2]{NR} we have
\[
e_0(\widehat M)=2e_0(M)\neq -1. 
\]
Thus, since $\widehat\Si_g=S^2$ and $\widehat M$ carries a positive,
transverse contact structure, we have $e_0(\widehat M)\leq
-2$. Therefore $e_0(M)\leq -1$, and case (a) of Theorem~\ref{t:main}
holds for $M$.

\sh{The subcase $g=0$, $e_0(M)=-1$, $r\geq 3$}

This is the last and hardest subcase. Before delving into the
proof, we need some preparation.

\sh{Algebraic preliminaries}

Let 
\[
\rho\in\Q_{>1}:=\{q\in\Q\ |\ q>1\}.
\]
Then, there is a unique way of writing $\rho$ as a continued fraction
\[
\rho=[a_1,\ldots,a_h],\quad a_i\geq 2,\ i=1,\ldots, h.
\]
It is well--known that if $\frac pq=[a_1,\ldots, a_h]$, with $p, q$
coprime, then $[a_h,\ldots,a_1]$ is of the form $\frac p{q'}$ for some
$q'$ coprime with $p$.

The set of finite sequences of integers has a natural linear order
$\preccurlyeq$ given by the following definition.

\begin{defi}
Declare 
\[
(a_1,\ldots, a_h)\preccurlyeq (b_1,\ldots,b_k) 
\]
if and only if there exists an index $1\leq j\leq\min(h,k)$ such that
\begin{itemize}
\item
$a_i=b_i$ for $i=1,\ldots, j-1$, and 
\item
either (i) $a_j < b_j$ or (ii) $j=k\leq h$ and $a_k=b_k$  .
\end{itemize}
\end{defi}

For example, we have
\[
(2,2,2)\preccurlyeq (2,2,3)\preccurlyeq (2,2)
\]

\begin{lem}\label{l:lexi=stand}
Let $\rho,\si\in\Q_{>1}$, with $\rho=[a_1,\ldots,a_h]$ and
$\si=[b_1,\ldots,b_k]$. Then,
\[
\rho\leq\si\quad\Longleftrightarrow\quad
(a_1,\ldots,a_h)\preccurlyeq (b_1,\ldots,b_k).
\]
\end{lem}

\begin{proof}
Suppose that $(a_1,\ldots,a_h)\preccurlyeq (b_1,\ldots,b_k)$.  Observe
that the rational expression $[a_1,\ldots,a_h]$ makes sense if
$a_1,\ldots,a_h\in\Q_{>1}$. Let $1\leq j\leq\min(h,k)$ be such that
$a_i=b_i$ for $1\leq i\leq j-1$ and either (i) $a_j < b_j$ or (ii)
$j=k\leq h$ and $a_k=b_k$. If (i) holds then
\[
[a_j,\ldots, a_h]\leq a_j \leq b_j-1 < [b_j,\ldots, b_k],
\]
therefore we have
\[
\begin{split}
\rho &= [a_1,\ldots, a_{j-1},[a_j,\ldots,a_h]] =\\ 
& = [b_1,\ldots, b_{j-1},[a_j,\ldots,a_h]]<
\si = [b_1,\ldots, b_{j-1},[b_j,\ldots, b_k]].
\end{split}
\]
If (ii) holds then, since $[b_k,a_{k+1},\ldots, a_h]\leq b_k$,
\[
\begin{split}
\rho &= [a_1,\ldots,a_k,a_{k+1},\ldots, a_h] = \\
&= [b_1,\ldots, b_{k-1},[b_k,a_{k+1},\ldots, a_h]] \leq 
\si = [b_1,\ldots, b_k].
\end{split}
\]

Conversely, suppose $\rho\leq\si$. If
$(a_1,\ldots,a_h)\not\preccurlyeq (b_1,\ldots,b_k)$ then
$(a_1,\ldots,a_h)\neq (b_1,\ldots,b_k)$, i.e. $\rho\neq\si$ and
$(b_1,\ldots,b_k)\preccurlyeq (a_1,\ldots,a_h)$, which implies $\si
<\rho$. Therefore we must have $(a_1,\ldots,a_h)\preccurlyeq
(b_1,\ldots,b_k)$.
\end{proof}

Consider the involution $\Q_{>1}\to\Q_{>1}$ which maps
$\rho\in\Q_{>1}$ to the only solution $\rho'\in\Q_{>1}$ of the
equation:
\[
\frac 1\rho + \frac 1{\rho'} = 1.
\]
Observe that $2=2'$, and $\rho_1\leq\rho_2$ if and only if
$\rho'_2\leq\rho'_1$. Therefore, the involution maps the interval
$(1,2)$ bijectively onto $(2,\infty)$ reversing the standard linear
order.

\begin{figure}[ht]
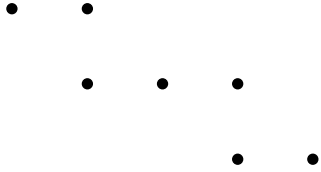
\small
\begin{center}
\begin{graph}(4,2)(0,2)
\graphnodesize{0.15}
  \roundnode{n1}(0,4)
  \roundnode{n2}(1,4)
  \roundnode{n3}(1,3)
  \roundnode{n4}(2,3)
  \roundnode{n5}(3,3)
  \roundnode{n6}(3,2)
  \roundnode{n7}(4,2)
\end{graph}
\end{center}
\caption{The point diagram $D([3,4,3],[2,3,2,3,2])$}
\label{f:fig2} 
\end{figure}

Given 
\[
\rho=[a_1,\ldots,a_h]\in\Q_{>1}, 
\]
the~\emph{Riemenschneider's point diagram}~\cite{Ri} $D(\rho,\rho')$
says how to compute the coefficients in the expansion
\[
\rho'=[b_1,\ldots,b_k]. 
\]
The diagram $D(\rho,\rho')$ consists of $h$ rows of dots, with the
$i$--th row consisting of $a_i-1$ dots and whose first dot lies under
the last dot of the $(i-1)$--st row. Then, $b_j$ is given by the
number of elements in the $j$--th column of $D(\rho,\rho')$ increased
by one. For example, if $\rho=[3,4,3]$ and $\rho'=[2,3,2,3,2]$, the
diagram is given by Figure~\ref{f:fig2}.

\sh{Setting up the stage for the proof}

\begin{figure}[ht]
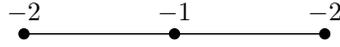
\small
\begin{center}
\begin{graph}(5,1)(-0.5,-0.2)
\graphnodesize{0.15}
  \roundnode{n1}(0,0)
  \roundnode{n2}(2,0)
  \roundnode{n3}(4,0)

  \edge{n1}{n2}
  \edge{n2}{n3}

  \autonodetext{n1}[n]{$-2$}
  \autonodetext{n2}[n]{$-1$}
  \autonodetext{n3}[n]{$-2$}

\end{graph}
\end{center}
\caption{Configuration of spheres in $X_M$ when $\de_2\leq 2$}
\label{f:fig2.5} 
\end{figure}

Recall that we are assuming $g=0$, $e_0(M)=-1$ and $r\geq 3$. We need
to show that if $M$ carries a positive, transverse contact structure
then $\Ga(M)=(\ga_1,\ldots,\ga_r)$ is realizable. Without loss of
generality, we may assume
\[
\ga_1\geq\ga_2\geq\cdots\geq\ga_r.
\]
Define 
\[
\de_i:=\frac 1{\ga_i},\ i=1,\ldots, r.
\]
The realizability of $\Ga(M)$ is equivalent to the existence of
coprime integers $a, m$ such that $1\leq a<m$ and
\begin{equation}\label{e:realizability}
\de_1>\frac ma,\quad \de_2>\frac
m{m-a},\quad \de_3,\ldots,\de_r > m.
\end{equation}
If $\de_1>2$ then Inequalities~\eqref{e:realizability} are satisfied
for $m=2$, $a=1$, and $\Ga(M)$ is realizable. Therefore we may assume
$\de_1\leq 2$. Under this assumption, we can write the continued
fraction expansion of $\de_1$ as:
\begin{equation}\label{e:del1}
\de_1=[\overbrace{2,\ldots,2}^{n_1+1}, n_2, \overbrace{2,\ldots,2}^{n_3},
n_4,\ldots, n_{2p}, \overbrace{2,\ldots,2}^{n_{2p+1}}]
\end{equation}
for some 
\[
n_1,n_3,\ldots,n_{2p+1}\geq 0,\quad
n_{2},n_4,\ldots,n_{2p}\geq 3.
\]
Next, we claim that we may also assume $\de_2>2$. In fact, if $\de_2\leq 2$
then 
\[
\de_2=[2,\ldots], 
\]
and inspecting Figure~\ref{f:fig0} we see that $X_M$ must contain a
configuration of spheres dual to the graph of Figure~\ref{f:fig2.5}.
By blowing down (-1)--spheres, this immediately implies that $X_M$
contains an embedded sphere of square zero, which is impossible by
Corollary~\ref{c:adjunction}. 

Assuming $\de_2>2$, we have
\begin{equation}\label{e:del2}
\de_2=[m_1,\overbrace{2,\ldots,2}^{m_2}, m_3,\ldots,
m_{2q-1},\overbrace{2,\ldots,2}^{m_{2q}}],
\end{equation}
for some 
\[
m_1,m_3,\ldots,m_{2q-1}\geq 3,\quad
m_2,m_4,\ldots,m_{2q}\geq 0.
\]
Since $\de_3\geq\de_2$, we also have
\[
\de_3=[d,\ldots],\quad\text{with}\quad d\geq 3.
\]

\sh{Blowing down and adjunction inequalities}

\begin{figure}[ht]
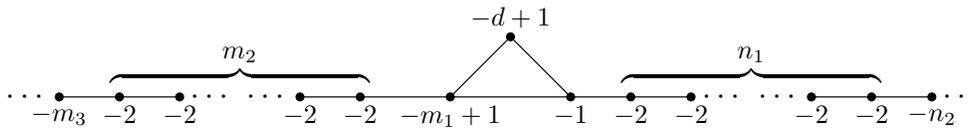
\small
\unitlength=.8cm
\begin{center}
\begin{graph}(15.5,2)(-1.5,-0.4)
\graphnodesize{0.15}
  \roundnode{m3}(-1,0)
  \roundnode{m21}(0,0)
  \roundnode{m22}(1,0)

  \roundnode{m23}(3,0)
  \roundnode{m24}(4,0)
  \roundnode{m1}(5.5,0)

  \roundnode{top}(6.5,1)

  \roundnode{n0}(7.5,0)
  \roundnode{n11}(8.5,0)
  \roundnode{n12}(9.5,0)

  \roundnode{n13}(11.5,0)
  \roundnode{n14}(12.5,0)
  \roundnode{n2}(13.5,0)

  \freetext(2,0.35){$\overbrace{\hspace{98pt}}$}
  \freetext(2,0.7){$m_2$}
  \freetext(10.5,0.35){$\overbrace{\hspace{98pt}}$}
  \freetext(10.5,0.7){$n_1$}

  \autonodetext{m3}[w]{\large $\cdots$}
  \edge{m3}{m21}
  \edge{m21}{m22}
  \autonodetext{m22}[e]{\large $\cdots$}

  \autonodetext{m23}[w]{\large $\cdots$}
  \edge{m23}{m24}
  \edge{m24}{m1}
  \edge{m1}{top}
  \edge{m1}{n0}
  \edge{top}{n0}
  \edge{n0}{n11}
  \edge{n11}{n12}
  \autonodetext{n12}[e]{\large $\cdots$}

  \autonodetext{n13}[w]{\large $\cdots$}
  \edge{n13}{n14}
  \edge{n14}{n2}
  \autonodetext{n2}[e]{\large $\cdots$}

  \autonodetext{m3}[s]{\small $-m_3$}
  \autonodetext{m21}[s]{\small $-2$}
  \autonodetext{m22}[s]{\small $-2$}
  \autonodetext{m23}[s]{\small $-2$}
  \autonodetext{m24}[s]{\small $-2$}
  \autonodetext{m1}[s]{\small $-m_1+1$}
  \autonodetext{top}[n]{\small $-d+1$}
  \autonodetext{n0}[s]{\small $-1$}
  \autonodetext{n11}[s]{\small $-2$}
  \autonodetext{n12}[s]{\small $-2$}
  \autonodetext{n13}[s]{\small $-2$}
  \autonodetext{n14}[s]{\small $-2$}
  \autonodetext{n2}[s]{\small $-n_2$}
\end{graph}
\end{center}
\caption{The manifold $X^{(0)}$}
\label{f:fig4} 
\end{figure}

\begin{figure}[ht]
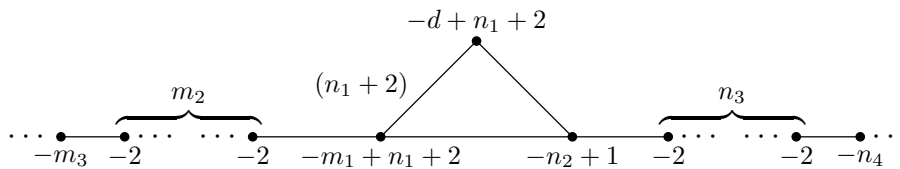
\small
\unitlength=0.85cm
\begin{center}
\begin{graph}(13,2.5)(-4,-0.4)
  \graphnodesize{0.15}

  \roundnode{m3}(-4,0)
  \roundnode{m21}(-3,0)
  \roundnode{m24}(-1,0)
  \roundnode{m1}(1,0)
  \roundnode{top}(2.5,1.5)
  \roundnode{n0}(4,0)
  \roundnode{n11}(5.5,0)
  \roundnode{n14}(7.5,0)
  \roundnode{n4}(8.5,0)

  \edge{m3}{m21}
  \edge{m24}{m1}
  \edge{m1}{top}
  \edge{m1}{n0}
  \edge{top}{n0}
  \edge{n0}{n11}
  \edge{n14}{n4}

  \autonodetext{m3}[w]{\large $\cdots$}
  \autonodetext{m21}[e]{\large $\cdots$}
  \autonodetext{m24}[w]{\large $\cdots$}
  \autonodetext{n11}[e]{\large $\cdots$}
  \autonodetext{n14}[w]{\large $\cdots$}
  \autonodetext{n4}[e]{\large $\cdots$}

  \freetext(-2,0.3){$\overbrace{\hspace{55pt}}$}
  \freetext(-2,0.65){\small $m_2$}
  \freetext(6.5,0.3){$\overbrace{\hspace{55pt}}$}
  \freetext(6.5,0.65){\small $n_3$}

  \autonodetext{m3}[s]{\small $-m_3$}
  \autonodetext{m21}[s]{\small $-2$}
  \autonodetext{m24}[s]{\small $-2$}
  \autonodetext{m1}[s]{\small $-m_1+n_1+2$}
  \autonodetext{top}[n]{\small $-d+n_1+2$}
  \autonodetext{n0}[s]{\small $-n_2+1$}
  \autonodetext{n11}[s]{\small $-2$}
  \autonodetext{n14}[s]{\small $-2$}
  \autonodetext{n4}[s]{\small $-n_4$}
  \freetext(0.7,0.8){\small $(n_1+2)$}

\end{graph}
\end{center}
\caption{The manifold $X^{(1)}$}
\label{f:fig5} 
\end{figure}

Since the central sphere in Figure~\ref{f:fig0} is a $(-1)$--sphere,
we can blow it down. This shows that the 4--manifold $X_M$ of
Section~\ref{s:filling} contains $X^{(0)}\#\overline\CP^2$, where
$X^{(0)}$ is the plumbing associated to the graph of
Figure~\ref{f:fig4}.
Blowing down $X^{(0)}$ $n_1+1$ times, we obtain Figure~\ref{f:fig5},
representing a four--manifold $X^{(1)}$ such that 
\[
X^{(0)}\cong X^{(1)}\# (n_1+1)\overline\CP^2. 
\]
Notice that each vertex in Figure~\ref{f:fig5} corresponds to an
embedded sphere. In this picture and the following ones, a numerical
weight of type ``$(n)$'' on the edge between two vertices of the graph
denotes the intersection number between suitably chosen homology
classes corresponding to the vertices.

By Corollary~\ref{c:adjunction}, $X_M$ contains no spheres of
non--negative self--intersection. Therefore, since $X^{(1)}$ can be
embedded in $X_M$, we have
\[
-m_1 + n_1 +2 \leq -1. 
\]

\begin{figure}[ht]
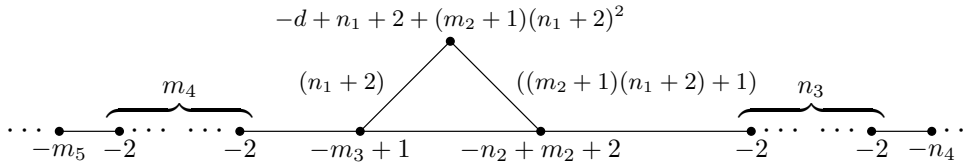
\small
\unitlength=0.8cm
\begin{center}
\begin{graph}(15.5,2.5)(-4.5,-0.4)
  \graphnodesize{0.15}

  \roundnode{m3}(-4,0)
  \roundnode{m21}(-3,0)
  \roundnode{m24}(-1,0)
  \roundnode{m1}(1,0)
  \roundnode{top}(2.5,1.5)
  \roundnode{n0}(4,0)
  \roundnode{n11}(7.5,0)
  \roundnode{n14}(9.5,0)
  \roundnode{n4}(10.5,0)

  \edge{m3}{m21}
  \edge{m24}{m1}
  \edge{m1}{top}
  \edge{m1}{n0}
  \edge{top}{n0}
  \edge{n0}{n11}
  \edge{n14}{n4}

  \autonodetext{m3}[w]{\large $\cdots$}
  \autonodetext{m21}[e]{\large $\cdots$}
  \autonodetext{m24}[w]{\large $\cdots$}
  \autonodetext{n11}[e]{\large $\cdots$}
  \autonodetext{n14}[w]{\large $\cdots$}
  \autonodetext{n4}[e]{\large $\cdots$}

  \freetext(-2,0.35){$\overbrace{\hspace{55pt}}$}
  \freetext(-2,0.7){\small $m_4$}
  \freetext(8.5,0.35){$\overbrace{\hspace{55pt}}$}
  \freetext(8.5,0.7){\small $n_3$}

  \autonodetext{m3}[s]{\small $-m_5$}
  \autonodetext{m21}[s]{\small $-2$}
  \autonodetext{m24}[s]{\small $-2$}
  \autonodetext{m1}[s]{\small $-m_3+1$}
  \autonodetext{top}[n]{\footnotesize $-d+n_1+2+(m_2+1)(n_1+2)^2$}
  \autonodetext{n0}[s]{\small  $-n_2+m_2+2$}
  \autonodetext{n11}[s]{\small $-2$}
  \autonodetext{n14}[s]{\small $-2$}
  \autonodetext{n4}[s]{\small $-n_4$}
  \freetext(0.7,0.8){\footnotesize  $(n_1+2)$}
  \freetext(5.6,0.8){\footnotesize  $((m_2+1)(n_1+2)+1)$}

\end{graph}
\end{center}
\caption{The manifold $X^{(2)}$}
\label{f:fig6} 
\end{figure}

\begin{figure}[ht]
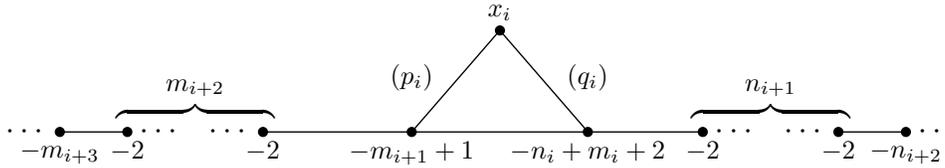
\small
\unitlength=0.9cm
\begin{center}
\begin{graph}(13,2.5)(-4.3,-0.4)
  \graphnodesize{0.15}

  \roundnode{m3}(-4,0)
  \roundnode{m21}(-3,0)
  \roundnode{m24}(-1,0)
  \roundnode{m1}(1.2,0)
  \roundnode{top}(2.5,1.5)
  \roundnode{n0}(3.8,0)
  \roundnode{n11}(5.5,0)
  \roundnode{n14}(7.5,0)
  \roundnode{n4}(8.5,0)

  \edge{m3}{m21}
  \edge{m24}{m1}
  \edge{m1}{top}
  \edge{m1}{n0}
  \edge{top}{n0}
  \edge{n0}{n11}
  \edge{n14}{n4}

  \autonodetext{m3}[w]{\large $\cdots$}
  \autonodetext{m21}[e]{\large $\cdots$}
  \autonodetext{m24}[w]{\large $\cdots$}
  \autonodetext{n11}[e]{\large $\cdots$}
  \autonodetext{n14}[w]{\large $\cdots$}
  \autonodetext{n4}[e]{\large $\cdots$}

  \freetext(-2,0.3){$\overbrace{\hspace{60pt}}$}
  \freetext(-2,0.65){\small $m_{i+2}$}
  \freetext(6.5,0.3){$\overbrace{\hspace{60pt}}$}
  \freetext(6.5,0.65){\small $n_{i+1}$}

  \autonodetext{m3}[s]{\small $-m_{i+3}$}
  \autonodetext{m21}[s]{\small $-2$}
  \autonodetext{m24}[s]{\small $-2$}
  \autonodetext{m1}[s]{\small $-m_{i+1}+1$}
  \autonodetext{top}[n]{\small $x_{i}$}
  \autonodetext{n0}[s]{\small $-n_i+m_i+2$}
  \autonodetext{n11}[s]{\small $-2$}
  \autonodetext{n14}[s]{\small $-2$}
  \autonodetext{n4}[s]{\small $-n_{i+2}$}
  \freetext(1.2,0.8){\small $(p_{i})$}
  \freetext(3.8,0.8){\small $(q_{i})$}
\end{graph}
\end{center}
\caption{The manifold $X^{(i)}$ for even $i$}
\label{f:fig7} 
\end{figure}

\begin{figure}[ht]
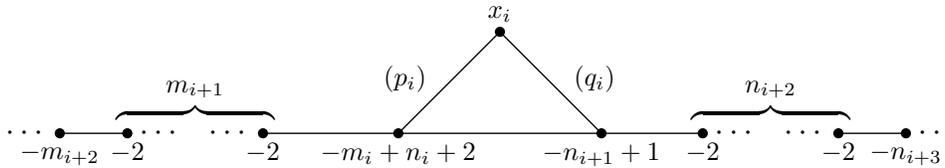
\small
\unitlength=0.9cm
\begin{center}
\begin{graph}(13,2.5)(-4.3,-0.4)
  \graphnodesize{0.15}

  \roundnode{m3}(-4,0)
  \roundnode{m21}(-3,0)
  \roundnode{m24}(-1,0)
  \roundnode{m1}(1,0)
  \roundnode{top}(2.5,1.5)
  \roundnode{n0}(4,0)
  \roundnode{n11}(5.5,0)
  \roundnode{n14}(7.5,0)
  \roundnode{n4}(8.5,0)

  \edge{m3}{m21}
  \edge{m24}{m1}
  \edge{m1}{top}
  \edge{m1}{n0}
  \edge{top}{n0}
  \edge{n0}{n11}
  \edge{n14}{n4}

  \autonodetext{m3}[w]{\large $\cdots$}
  \autonodetext{m21}[e]{\large $\cdots$}
  \autonodetext{m24}[w]{\large $\cdots$}
  \autonodetext{n11}[e]{\large $\cdots$}
  \autonodetext{n14}[w]{\large $\cdots$}
  \autonodetext{n4}[e]{\large $\cdots$}

  \freetext(-2,0.3){$\overbrace{\hspace{60pt}}$}
  \freetext(-2,0.65){\small $m_{i+1}$}
  \freetext(6.5,0.3){$\overbrace{\hspace{60pt}}$}
  \freetext(6.5,0.65){\small $n_{i+2}$}

  \autonodetext{m3}[s]{\small $-m_{i+2}$}
  \autonodetext{m21}[s]{\small $-2$}
  \autonodetext{m24}[s]{\small $-2$}
  \autonodetext{m1}[s]{\small $-m_i+n_i+2$}
  \autonodetext{top}[n]{\small $x_{i}$}
  \autonodetext{n0}[s]{\small $-n_{i+1}+1$}
  \autonodetext{n11}[s]{\small $-2$}
  \autonodetext{n14}[s]{\small $-2$}
  \autonodetext{n4}[s]{\small $-n_{i+3}$}
  \freetext(1.1,0.8){\small $(p_{i})$}
  \freetext(3.9,0.8){\small $(q_{i})$}
\end{graph}
\end{center}
\caption{The manifold $X^{(i)}$ for odd $i$}
\label{f:fig8} 
\end{figure}

If $-m_1+n_1+2 = -1$ we blow down another $m_2+1$ times obtaining 
a 4--manifold $X^{(2)}$ such that 
\[
X^{(1)}\cong X^{(2)}\# (m_2+1)\overline\CP^2, 
\]
with $X^{(2)}$ given by Figure~\ref{f:fig6}. In Figure~\ref{f:fig6}
the vertices on the horizontal line correspond to embedded spheres and
the top vertex can be represented by an immersed sphere with
$(m_2+1)\binom{n_1+2}{2}$ positive, transverse
self--intersections. Observe that, by smoothing out the
self--intersections, we can represent the same homology class by an
embedded surface of genus
\[
g_2=(m_2+1)\binom{n_1+2}{2}. 
\]
Moreover, by Corollary~\ref{c:adjunction} we have
\[
-n_2+m_2+2\leq -1.
\]
Continuing in the same fashion, this process gives rise
to a sequence of manifolds
\[
X^{(0)}\to X^{(1)}\to X^{(2)}\to\cdots\to X^{(i)}\to\cdots 
\]
such that each $X^{(i)}$ (if defined) is given by Figure~\ref{f:fig7}
for even $i$ and by Figure~\ref{f:fig8} for odd $i$. By
Corollary~\ref{c:adjunction}, we have
\[
(-1)^{i}(m_i -n_i) + 2 \leq -1
\]
for each $i$ such that $X^{(i)}$ is defined. Moreover, if 
$X^{(i)}$ is defined, then in order for $X^{(i+1)}$ to be 
defined as well we need to have 
\begin{equation}\label{e:m-n}
(-1)^{i}(m_i -n_i) + 2 = -1.
\end{equation} 
If Equation~\eqref{e:m-n} holds, then 
\[
X^{(i)}\cong
\begin{cases}
X^{(i+1)}\#(n_{i+1}+1)\overline\CP^2\quad\text{if $i$ is even},\\
X^{(i+1)}\#(m_{i+1}+1)\overline\CP^2\quad\text{if $i$ is odd}.\\
\end{cases}
\]
As in the case of $X^{(2)}$, the homology class corresponding to the
top vertex in each of Figures~\ref{f:fig7} and~\ref{f:fig8} can be
represented in a natural way by an immersed sphere $S_i\subset
X^{(i)}$. We denote by $g_i$ the genus of the smooth surface $\Si_i$
obtained by smoothing out the singularities of $S_i$. Then, the
numbers $p_i$ and $q_i$ denote algebraic as well as geometric
intersection numbers between $\Si_i$ and two embedded spheres
representing the homology classes which correspond to the vertices
connected to the top. Moreover, it is easy to check that the
following relations hold:
\begin{equation}\label{e:even-odd}
\begin{split}
x_{i+1} &=
\begin{cases}
x_i+(n_{i+1}+1)q_{i}^2\quad\text{if $i$ is even}\\
x_i+(m_{i+1}+1) p_{i}^2\quad\text{if $i$ is odd}
\end{cases}\\
p_{i+1} &=
\begin{cases}
p_i+(n_{i+1}+1) q_{i}\quad\text{if $i$ is even}\\
p_{i}\quad\text{if $i$ is odd}
\end{cases}\\
q_{i+1} &= 
\begin{cases}
q_{i}\quad\text{if $i$ is even}\\
q_{i} + (m_{i+1}+1)p_{i}\quad\text{if $i$ is odd}
\end{cases}\\
g_{i+1} &=
\begin{cases}
g_i+(n_{i+1} + 1)\binom{q_{i}}{2}\quad\text{if $i$ is even}\\
g_i+(m_{i+1} + 1)\binom{p_{i}}{2}\quad\text{if $i$ is odd}
\end{cases}
\end{split}
\end{equation}
Strictly speaking, Figures~\ref{f:fig7} and~\ref{f:fig8} should be
taken literally only if 
\begin{equation}\label{e:range}
0\leq i<\min(2q,2p+1). 
\end{equation}
(Recall that $p$ and $q$ were defined in Equations~\eqref{e:del1}
and~\eqref{e:del2}, respectively). When $i=\min(2q,2p+1)$, the
pictures should be suitably interpreted, because the diagram has only
one horizontal leg. Similarly, Equations~\eqref{e:even-odd} shold be
thought of as relations only when~\eqref{e:range} holds, while for
$i=\min(2q,2p+1)$ they should be thought of as the definition of
$x_{i+1}$, $p_{i+1}$, $q_{i+1}$ and $g_{i+1}$.

Using~\eqref{e:even-odd} it is easy to check that the number
\[
2g_i - 2 - x_i + p_i + q_i,\quad 0\leq i\leq\min(2q+1,2p+2), 
\]
is independent of $i$, and therefore it is always equal to its value
for $i=0$, i.e.  $d-1$. On the other hand, the adjunction inequality
implies
\[
2g_i -2 -x_i\geq 0, 
\]
hence
\begin{equation}\label{e:d}
d > p_i + q_i,\quad 0\leq i\leq\min(2q+1,2p+2).
\end{equation}

\sh{The end of the proof}

By considering the Riemenschneider's point diagram $D(\de_2,\de_2')$,
we see that
\[
\de_2'=[\overbrace{2,\ldots,2}^{m_1-2}, m_2+3, \overbrace{2,\ldots,2}^{m_3-3},
m_4+3,\ldots,\overbrace{2,\ldots,2}^{m_{2q-1}-3},m_{2q}+2].
\]

Clearly, one of the following holds: 
\begin{enumerate}
\item
There exists $k$ such that $0\leq k\leq\min(2q,2p+1)$,
Equation~\eqref{e:m-n} holds for every $0\leq i<k$, and 
\[
(-1)^k(m_k -n_k) + 2 < -1.
\]
\item
Equation~\eqref{e:m-n} holds for every $0\leq i\leq\min(2q,2p+1)$.
\end{enumerate}

We shall now treat separately the possible cases which can occur.

\sh{(1) holds and $k$ is even}

In this case we have
\begin{multline*}
\de_2' <
[\overbrace{2,\ldots,2}^{m_1-2},m_2+3,\overbrace{2,\ldots,2}^{m_3-3},
\ldots, m_{k-2} + 3,\overbrace{2,\ldots,2}^{m_{k-1}-3}, m_k+3] \leq
\\ \leq [\overbrace{2,\ldots,2}^{n_1+1},n_2,\overbrace{2,\ldots,2}^{n_3},
\ldots, n_{k-2},\overbrace{2,\ldots,2}^{n_{k-1}}, n_i-1] < \\
[\overbrace{2,\ldots,2}^{n_1+1},n_2,\overbrace{2,\ldots,2}^{n_3},
\ldots, n_{k-2},\overbrace{2,\ldots,2}^{n_{k-1}}, n_i,\ldots ]
=\de_1
\end{multline*}

\begin{lem}\label{l:cal}
Let $0\leq k\leq\min(2q+1,2p+2)$, and suppose that Equation~\eqref{e:m-n}
holds for every $0\leq i<k$. Then, if $k$ is even we have
\begin{equation*}
\frac{p_k+q_k}{p_k}=
[m_k+3,\overbrace{2,\ldots, 2}^{m_{k-1}-3},m_{k-2}+3, \ldots,
\overbrace{2,\ldots,2}^{m_3 - 3}, m_2 + 3,\overbrace{2,\ldots,2}^{m_1 - 2}],
\end{equation*}
while if $k$ is odd then
\begin{equation*}
\frac{p_k+q_k}{q_k} = 
[\overbrace{2,\ldots,2}^{n_k+1}, n_{k-1},
\overbrace{2,\ldots,2}^{n_{k-2}}, n_{k-3}, \ldots, 
\overbrace{2,\ldots,2}^{n_3}, n_2,\overbrace{2,\ldots,2}^{n_1+1}].
\end{equation*}
\end{lem}

\proof
Observe that Relations~\eqref{e:even-odd} imply
\begin{equation}\label{e:fractions}
\frac{p_{i+1}}{q_{i+1}}=\frac{p_{i}}{q_{i}} + 
\begin{cases}
(n_{i+1} + 1)\quad\text{if $i$ is even},\\
(m_{i+1} + 1)\quad\text{if $i$ is odd},
\end{cases}\quad
0\leq i<k.
\end{equation}
Let us introduce the notation:
\[
[n_1,\ldots, n_k]^+ := 
n_1 + \cfrac{1}{n_2 +
          \cfrac{1}{\ddots +
             \cfrac{1}{n_k
}}}
\]
Since Equation~\eqref{e:m-n} holds for every $0\leq i<k$,
by~\eqref{e:fractions} we get, when $k$ is even,
\begin{equation*}
1+\frac{q_k}{p_k} = 
[ m_k+2, m_{k-1}-2, m_{k-2}+1,\ldots, m_3 - 2, m_2+1, m_1-1]^+,
\end{equation*}
and when $i$ is odd 
\begin{equation*}
1+\frac{p_k}{q_k} = 
[n_k+2, n_{k-1}-2, n_{k-2}+1,\ldots, n_3 +1, n_2-1, n_1+2]^+.
\end{equation*}
The lemma follows from the following identities, which can be
established by a straightforward induction:
\[
\begin{split}
[m_k+2, m_{k-1}-2, m_{k-2}+1,\ldots, m_3 - 2, m_2 + 1, m_1-1]^+ =\\
[m_k+3,\overbrace{2,\ldots, 2}^{m_{k-1}-3},m_{k-2}+3, \ldots,
\overbrace{2,\ldots,2}^{m_3 - 3}, m_2 + 3,\overbrace{2,\ldots,2}^{m_1 - 2}].
\end{split}
\]
$$[ n_{k+1}+2, n_k-2, n_{k-1}+1,\ldots, n_3 +1, n_2-1, n_1+2]^+ =$$
$$[\overbrace{2,\ldots,2}^{n_{k+1}+1}, n_k,
\overbrace{2,\ldots,2}^{n_{k-1}}, n_{k-2}, \ldots, 
\overbrace{2,\ldots,2}^{n_3}, n_2,\overbrace{2,\ldots,2}^{n_1+1}].\eqno{\qed}$$
\eject

By~\eqref{e:d} and Lemma~\ref{l:cal}, $\de_3$ is bigger than the
numerator of a fraction representing the number
\[
\rho =
[\overbrace{2,\ldots,2}^{m_1-2},m_2+3,\overbrace{2,\ldots,2}^{m_3-3},
\ldots, m_{k-2} + 3,\overbrace{2,\ldots,2}^{m_{k-1}-3}, m_k+3].
\]
If $\rho=\frac ma$ with $a$ and $m$ coprime, then 
\[
\de_1>\frac ma,\quad \de_2>\rho'=\frac m{m-a}
\quad\text{and}\quad \de_3>m. 
\]
Therefore $\Ga(M)$ is realizable.

\sh{(1) holds and $k$ is odd}

In this case we have
\begin{multline*}
\de_2' <
[\overbrace{2,\ldots,2}^{m_1-2},m_2+3,\overbrace{2,\ldots,2}^{m_3-3},
\ldots, m_{k-3} + 3,\overbrace{2,\ldots,2}^{m_{k-2}-3}, m_{k-1}+3,
\overbrace{2,\ldots,2}^{m_k-3}] \leq
\\ \leq [\overbrace{2,\ldots,2}^{n_1+1},n_2,\overbrace{2,\ldots,2}^{n_3},
\ldots, n_{k-3},\overbrace{2,\ldots,2}^{n_{k-2}}, n_{k-1},
\overbrace{2,\ldots,2}^{n_k+1}] < \\
[\overbrace{2,\ldots,2}^{n_1+1},n_2,\overbrace{2,\ldots,2}^{n_3},
\ldots, n_{k-3},\overbrace{2,\ldots,2}^{n_{k-2}}, n_{k-1},
\overbrace{2,\ldots,2}^{n_k},\ldots ]
=\de_1
\end{multline*}
By~\eqref{e:d} and Lemma~\ref{l:cal}, $\de_3$ is bigger than the
numerator of
\[
\frac ma =
[\overbrace{2,\ldots,2}^{n_1+1},n_2,\overbrace{2,\ldots,2}^{n_3},
\ldots, n_{k-3},\overbrace{2,\ldots,2}^{n_{k-2}}, n_{k-1},
\overbrace{2,\ldots,2}^{n_k+1}].
\]
As in the previous case, it follows that $\Ga(M)$ is realizable.

\sh{(2) holds and $2q<2p+1$}

We have 
\begin{multline*}
\de_2' <
[\overbrace{2,\ldots,2}^{m_1-2},m_2+3,
\ldots, m_{2q}+3,\overbrace{2,\ldots,2}^{n_{2q+1}+1}] =
\\ = [\overbrace{2,\ldots,2}^{n_1+1},n_2,
\ldots, n_{2q},\overbrace{2,\ldots,2}^{n_{2q+1}+1}] < \\
< [\overbrace{2,\ldots,2}^{n_1+1},n_2,
\ldots, n_{2q},\overbrace{2,\ldots,2}^{n_{2q+1}},\ldots ]
=\de_1.
\end{multline*}
By~\eqref{e:d} and Lemma~\ref{l:cal}, $\de_3$ is bigger than the 
numerator of 
\[
[\overbrace{2,\ldots,2}^{n_1+1},n_2,
\ldots, n_{2q},\overbrace{2,\ldots,2}^{n_{2q+1}+1}],
\]
therefore $\Ga(M)$ is realizable.

\sh{(2) holds and $2p+1<2q$}

We have
\begin{multline*}
\de_2' <
[\overbrace{2,\ldots,2}^{m_1-2},m_2+3,
\ldots, m_{2p}+3,\overbrace{2,\ldots,2}^{m_{2p+1}-3}, m_{2p+2}+3] =
\\ = [\overbrace{2,\ldots,2}^{n_1+1},n_2,
\ldots, n_{2p},\overbrace{2,\ldots,2}^{n_{2p+1}}, m_{2p+2}+3] < \\
< [\overbrace{2,\ldots,2}^{n_1+1},n_2,
\ldots, n_{2p},\overbrace{2,\ldots,2}^{n_{2p+1}}]=\de_1.
\end{multline*}
By~\eqref{e:d} and Lemma~\ref{l:cal}, $\de_3$ is bigger than the 
numerator of 
\[
[\overbrace{2,\ldots,2}^{m_1-2},m_2+3,
\ldots, m_{2p}+3,\overbrace{2,\ldots,2}^{m_{2p+1}-3}, m_{2p+2}+3],
\]
therefore $\Ga(M)$ is realizable.

\Addresses\recd

\end{document}